\documentclass[aoas,preprint]{imsart}

\RequirePackage[OT1]{fontenc}
\RequirePackage{amsthm,amsmath}
\RequirePackage[numbers]{natbib}
\RequirePackage[colorlinks,citecolor=blue,urlcolor=blue]{hyperref}
\usepackage[pdftex]{graphicx}

\startlocaldefs
\numberwithin{equation}{section}
\theoremstyle{plain}
\newtheorem{Theorem}{Theorem}[section]
\newtheorem{Remark}{Remark}[section]
\newtheorem{Lemma}{Lemma}[section]
\newtheorem{Corollary}{Corollary}[section]
\endlocaldefs

\begin{document}

\begin{frontmatter}
\title{Gaussian One-Armed Bandit and Optimization of Batch Data
Processing} \runtitle{Gaussian One-Armed Bandit}

\begin{aug}
\author{\fnms{Alexander} \snm{Kolnogorov}\thanksref{m1}\ead[label=e1]{Alexander.Kolnogorov@novsu.ru}},

\runauthor{A. Kolnogorov}

\affiliation{Yaroslav-the-Wise Novgorod State
University\thanksmark{m1}}

\address{41 B.Saint-Petersburgskaya Str.,
Velikiy Novgorod, Rusiia, 173003\\
Applied Mathematics and Information Science Department\\
\printead{e1}}

\end{aug}

\begin{abstract}
We consider the minimax setup for Gaussian one-armed bandit
problem, i.e. the two-armed bandit problem with Gaussian
distributions of incomes and known distribution corresponding to
the first arm. This setup naturally arises when the optimization
of batch data processing is considered and there are two
alternative processing methods available with a priori known
efficiency of the first method. One should estimate the efficiency
of the second method and provide predominant usage of the most
efficient of both them. According to the main theorem of the
theory of games minimax strategy and minimax risk are searched for
as Bayesian ones corresponding to the worst-case prior
distribution.

As a result, we obtain the recursive integro-difference equation
and the second order partial differential equation in the limiting
case as the number of batches goes to infinity. This makes it
possible to determine minimax risk and minimax strategy by
numerical methods. If the number of batches is large enough we
show that batch data processing almost does not influence the
control performance, i.e. the value of the minimax risk. Moreover,
in case of Bernoulli incomes and large number of batches, batch
data processing provides almost the same minimax risk as the
optimal one-by-one data processing.
\end{abstract}



\end{frontmatter}

\def \eps{\varepsilon}
\def \mE{\mathrm{E}}
\def \tR{\tilde{R}}
\def \tg{\tilde{g}}
\def \hD {\hat{D}}
\def \hX {\hat{X}}
\def \hM {\hat{M}}

\def \hm {\hat{m}}
\def \hlambda {\hat{\lambda}}
\def \hTheta {\hat{\Theta}}
\def \hR {\hat{R}}
\def \hv {\hat{v}}
\def \hrho {\hat{\rho}}
\def \hL {\hat{L}}
\def \hC {\hat{C}}
\def \hsigma {\hat{\sigma}}
\def \tm {\tilde{m}}
\def \tR {\tilde{R}}
\def \tr {\tilde{r}}
\def \tg {\tilde{g}}
\def \tTheta {\tilde{\Theta}}

\section{Introduction}

The two-armed bandit problem originates from the slot machine with
two arms the choice of each is followed by the random income of
the gambler depending only on chosen arm. The goal is to maximize
the total expected income. To this end, the gambler should
determine more profitable arm and provide its predominant usage.
The problem has numerous applications in medical trials,
biological modelling, data processing, internet, etc.
(see~\cite{BF, PS} and references therein). It is also well-known
as the problem of expedient behavior in a random
environment~\cite{Tsetlin} and the problem of adaptive
control~\cite{Sragovich}.

In this article, we consider Gaussian one-armed bandit problem,
i.e., Gaussian two-armed bandit problem with known distribution of
income corresponding to the first action. Formally, let's consider
a controlled random process $\xi_n$, $n=1,\ldots,N$, which values
are interpreted as incomes, depend only on currently chosen
actions $y_n$ ($y_n\in\{1, 2\}$) and have Gaussian (normal)
distribution with a density $f_D(x| m)=(2\pi D)^{-1/2} \exp
\left\{-(x-m)^2/(2D)\right\}$ if $y_n=2$, where $D$ is assumed to
be known and $m$ is assumed to be unknown (later on, the
assumption of known $D$ can be omitted). If $y_n=1$ then
mathematical expectation $m_1$ is assumed to be known and without
loss of generality $m_1=0$ (otherwise, one can consider the
process $\xi_n-m_1$, $n=1,\ldots,N$).

A control strategy $\sigma$ at the point of time $k+n+1$ assigns,
in general, a random selection of actions depending on the current
history which is described by triplet $(k,X,n)$, where $k,n$ are
current cumulative numbers of the first and the second actions
applications, $X$ is current cumulative income for the use of the
second action. The value of cumulative income for the use of the
first action is immaterial because corresponding distribution is
known.

Considered random process is completely described by a vector
parameter~$\theta=(0,m)$. If a parameter $\theta$ was known then
the optimal strategy should always prescribe to choose the action
corresponding to the larger value of $0,m$. In this case, total
expected income would be equal to $N (0 \vee m)$ where $0 \vee m$
is the maximum of 0, $m$. And if the parameter is unknown then the
loss function
\begin{gather}\label{a1}
L_N(\sigma,\theta)=N (0 \vee
m)-\mE_{\sigma,\theta}\Biggl(\,\sum_{n=1}^N \xi_n\Biggr)
\end{gather}
describes losses of the total expected income with respect to its
maximal possible value because of the incomplete information. Here
$\mE_{\sigma,\theta}$ denotes the mathematical expectation with
respect to measure generated by strategy $\sigma$ and parameter
$\theta$. We assume that the set of parameters is
$\Theta=\{\theta=(0, m):\:|m|\le C\}$, where $0<C< \infty$ ensures
the boundedness of the loss function on $\Theta$. Definition
\eqref{a1} implies that one-armed bandit problem may be considered
as a game. In this game, the first player is the gambler with the
set of strategies $\{\sigma\}$. And the second player is the {\em
nature} with the set of strategies $\Theta$.

By the loss function \eqref{a1} we define the minimax risk
\begin{equation}\label{a2}
R^M_N(\Theta)=\inf_{\{\sigma\}} \sup_{\Theta} L_N(\sigma,\theta),
\end{equation}
the corresponding optimal strategy $\sigma^M$ is called the
minimax strategy.

Minimax approach to considered problem for Bernoulli two-armed
bandit was proposed in~\cite{Robbins}. It was shown in~\cite{FZ}
that explicit determination of the minimax strategy and minimax
risk for Bernoulli two-armed bandit is virtually impossible
already for $N>4$. However, in~\cite{Vogel} an asymptotic (as
$N\to\infty$) unimprovable in order bounds of the minimax risk
were obtained:
\begin{gather}\label{a3}
0.375 \le (DN)^{-1/2} R_N^M(\Theta)\le 0.752.
\end{gather}
Here $D=0.25$ is the maximum value of the variance of on-step
Bernoulli income. The bounds \eqref{a3} hold true for Gaussian
two-armed bandit as well. Obviously, the upper bound \eqref{a3}
remains valid for the one-armed bandit, too. Another approach to
robust control in the multi-armed bandit problem was considered
in~\cite{Nazin}. In this article mirror descent algorithm is used.
Note that in~\cite{Nazin} the minimization of total expected
income was considered. Therefore, instead of ``incomes'' they
considered ``losses'' and the loss function was chosen the
following
\begin{gather*}
L_N(\sigma,\theta)=\mE_{\sigma,\theta}\Biggl(\,\sum_{n=1}^N
\xi_n\Biggr)-N (0 \wedge m),
\end{gather*}
where $0 \wedge m$ is the minimum of 0, $m$. In case of Gaussian
distributions of incomes this setup can be reduced to presented in
this article by considering incomes $\{-\xi_n\}$,

Let's explain the choice of Gaussian distribution for incomes. We
consider the problem as applied to control of processing of large
amounts of data in sufficiently small number of stages by
partitioning them into batches. Let $T=NM$ items of data be given
which can be processed by one of two alternative methods and let
$\zeta_t$ denote the result of processing of the data item
numbered $t$. For example, processing can be successful
($\zeta_t=1$) or unsuccessful ($\zeta_t=0$) and one has to
maximize the total expected number of successfully processed data.
Or $\zeta_t$ is a duration of processing of the $t$-th data item
and one has to minimize the total expected computer time of data
processing. Assume that distributions of $\{\zeta_t\}$ depend only
on chosen methods and mathematical expectation $\mE(
\zeta_t|y_t=1)=m_1$ is known. Let's partition all the data into
$N$ batches, each containing $M$ items of data, and use the same
method for data processing in the same batch. For the control,
let's use the values of the process
$\xi_n=M^{-1/2}\sum\limits_{t=(n-1)M+1}^{nK} (\zeta_t-m_1)$,
$n=1,\ldots,N$. According to the central limit theorem
distributions of the process $\xi_n$, $n=1,\dots,N$ are close to
Gaussian with zero mathematical expectation corresponding to the
first method and variances are equal to one-step variances of
incomes just like in considered setting. In what follows, the
values of the process $\{\xi_n\}$ will be also considered as the
data, e.g., passed the preprocessing.

Note that the data in the same batches can be often processed in
parallel. In this case the total processing time depends on the
number of batches rather than on the total number of data.
However, there is a question of losses in control performance due
to such clustering of data. Numerical results given in
Section~\ref{num} show that the scaled minimax risk
$N^{-1/2}R^M_N(\Theta)$ is almost constant for comparatively small
$N$, e.g. for $N\ge 50$. Therefore, say, $50000$ items of data may
be processed with approximately equal maximal losses either in
1000 stages by batches of 50 data or in 50 stages by batches of
1000 data. However, it is to the full true only in the case of
close mathematical expectations $m_1$, $m_2$. If one is not sure
in closeness of $m_1$, $m_2$ then  sizes of batches at the initial
stage should be chosen smaller. Corresponding example is also
given in Section~\ref{num}.

\begin{Remark}
Parallel control in the two-armed bandit problem was first
proposed for the treatment a large group of patients by two
alternative drugs with different unknown efficiencies. Really, if
a doctor treats, say, one thousand patients one-by-one and the
result of treatment manifests in a week then the total treatment
would take about twenty years. Therefore, it was proposed to give
initially both drugs to sufficiently large test groups of patients
and then the more efficient drug to all the rest ones. As a
result, the treatment would take two weeks! The problem is
discussed in~\cite{LLRS} (see the references therein).
\end{Remark}

A very popular approach to the problem is a Bayesian one. Let's
denote by $\lambda(\theta)$ a prior probability distribution
density  on $\Theta$. The value
\begin{equation}\label{a4}
R^B_N(\lambda)= \min_{\{\sigma\}} \int_\Theta L_N(\sigma,\theta)
\,\lambda(\theta) d\theta
\end{equation}
is called the Bayesian risk and corresponding optimal strategy
$\sigma^B$ is called the Bayesian strategy. Bayesian approach
allows to find Bayesian risk and Bayesian strategy by dynamic
programming technique for any prior distribution. Both Bayesian
and minimax approaches are integrated by the main theorem of
theory of games according to which
\begin{equation}\label{a5}
R_N^M(\Theta)= R^B_N(\lambda^0)=\max_{\{\lambda\}} R^B_N(\lambda),
\end{equation}
i.e. minimax risk~\eqref{a2} is equal to Bayesian risk~\eqref{a4}
calculated with respect to the worst-case prior distribution
corresponding to the maximum of Bayesian risk. And minimax
strategy is equal to corresponding the Bayesian one. This theorem
for considered problem was proved in~\cite{Koln11} for even more
general setting.

The one-armed bandit problem with Bernoulli incomes was considered
in Bayesian setting in~\cite{BJK, Cher}. The main feature of the
Bayesian strategy which was proved in~\cite{BJK, Cher} is based on
the following idea. Since choosing the first action does not
influence the available information then, being once chosen, it
will be applied till the end of the control. Another proved
in~\cite{BJK, Cher} important feature of the strategy is its
thresholding property. We show below that in considered setting
these properties take place as well.

The structure of the article is the following. In
Section~\ref{equ} the recurrent integro-difference equations are
presented which allow to calculate Bayesian risk, Bayesian
strategy and expected losses for batch data processing. We prove
that batch data processing almost does not increase the minimax
risk if the number of batches is large enough. In
Section~\ref{lim} we present invariant notations of these
equations with a control horizon equal to unit. We prove the
existence and some properties of the limiting solution to
invariant equation and then present its description by the second
order partial differential equation. Usage of numerical methods in
Section~\ref{num} gives the following asymptotic estimate
\begin{gather}\label{a6}
\lim_{N \to \infty}(DN)^{-1/2} R_N^M(\Theta) =r
\end{gather}
with $r\approx 0.37$. Results of Section~\ref{num} make it
possible to omit the assumption of known $D$ if the number of data
items is large enough. First, the expected losses deviate just a
little if the variance $D$ is assigned with a significant error up
to 5\%. This means that unknown variance can be estimated at the
start of the control and then obtained estimate can be used.
Second, it turned out that maximal expected losses corresponding
to $D:0.25D_0\le D \le D_0$ are almost not more than those
corresponding to $D_0$, i.e. the minimax strategy corresponding to
$D_0$ saves approximately this property for all $D:0.25D_0\le D
\le D_0$.

Recall that batch data processing almost does not increase the
minimax risk if the number of batches is large enough. According
to the central limit theorem it means that the usage of batch data
processing makes it possible to ensure the value of the minimax
risk close to \eqref{a6} for a wide class of processes with the
same mathematical expectations and variances of one-step incomes.
However, this does not straightforwardly ensures  that minimax
risk cannot be diminished by the usage of one-by-one data
processing. In Section~\ref{asymp}, we show that one-by-one data
processing does not allow to diminish the minimax risk in
Bernoulli case. To calculate Bayesian risk in Bernoulli case, we
obtain the same partial differential equation as for Gaussian
one-armed bandit. Seemingly, this approach can be used for other
distributions of incomes, too.

Discussion is presented in Section~\ref{disc}.

Note that some results for Gaussian one-armed bandit corresponding
to $D=1$ were obtained in~\cite{Koln15}. Gaussian two-armed bandit
with different unknown variances $D_1$, $D_2$ was considered
in~\cite{Koln18}.

\section{Recursive Equations Describing Batch Data Processing Optimization}\label{equ}

According to the equality~\eqref{a5} we search minimax risk as
Bayesian one corresponding to the worst-case prior distribution.
In what follows, let's use control strategies which can change
actions after their application $M$ times in succession only. If
incomes of one-armed bandit $\{\xi_n\}$ come in succession then
such a strategy allows to switch actions more rarely. And if
incomes come in batches then this strategy allows the parallel
processing. For convenience we assume that $N$ is multiple of $M$.

Let's denote by $\lambda(m)$ a prior probability distribution
density of the second part of parameter $\theta=(0,m)$. We assume
that
\begin{equation}\label{b1}
\int_{-C}^0 |m|\,\lambda(m) dm>0,\qquad \int_0^C m\,\lambda(m)dm
>0.
\end{equation}
Let's suppose that the first and the second actions were chosen
$k$ and $n$ times. Then the control history up to the point of
time $k+n+1$ is described by statistics $(k,X,n)$ where $X$ is
cumulative income for the application of the second action
(cumulative income for the application of the first action does
not influence the available information and is not used). Denote
$n^*=nD$, $M^*=MD$. The posterior distribution density is then
equal to
\begin{gather}\label{b2}
\begin{aligned}
\lambda(m| X,n)&=\frac{f_{n^*}(X|n m)\,\lambda(m)} {P(X,n)}, \\
\mbox{where}\ &P(X,n)= \displaystyle{\int_{-C}^{C}} f_{n^*}(X|n
m)\,\lambda(m) dm.
\end{aligned}
\end{gather}
If additionally it is assumed that $X=0$ and $f_{n^*}(X|nm)=1$ as
$n=0$ then $\lambda(m| 0,0)=\lambda(m)$. Let's denote by
$R^B(k,X,n)$ Bayesian risk on the last $(N-k-n)$ steps calculated
with respect to posterior distribution density $\lambda(m| X,n)$.
The Bayesian risk~\eqref{a4} is then equal to $R^B(0,0,0)$ and can
be found by solving the standard recursive equation
\begin{equation}\label{b3}
R^B(k,X,n)=\min\bigl(R^B_1(k,X,n),R^B_2(k,X,n)\bigr),
\end{equation}
where $R^B_1(k,X,n)=R^B_2(k,X,n)=0$ for $k+n=N$ and
\begin{gather}\label{b4}
\begin{aligned}
R^B_1(k,X,n)&= \displaystyle{\int_0^C} M m \,\lambda(m| X,n)dm+R^B(k+M,X,n),\\
R^B_2(k,X,n)&= \displaystyle{\int_{-C}^0} \Big( |Mm|  \\
  &+\mE_Y R^B(k,X+Y,n+M)\Big) \lambda(m | X,n)dm
\end{aligned}
\end{gather}
for $k+n<N$. Here $\mE_Y$ denotes mathematical expectation with
respect to the density $f_{M^*}(Y|Mm)$:
$$
\mE_Y R(Y)= \int_{-\infty}^{+\infty} R(Y) f_{M^*}(Y|Mm) \,dY.
$$
In equations above, $R^B_\ell(\cdot)$ means cumulative expected
income on the residual control horizon of the length $(N-k-n)$ if
at first $M$ times the $\ell$-th action was chosen and then the
control was optimally implemented ($\ell=1,2$). Bayesian strategy
prescribes always to choose the action corresponding to the
smaller value of $R^B_1(\cdot),R^B_2(\cdot)$; the choice may be
arbitrary in case of their equality.

Consider a strategy
$\sigma:\sigma_\ell(k,X,n)=\Pr(y_{k+n+1}=\ell|k,X,n)$. Then
standard recursive equation for losses  takes the form
\begin{equation}\label{b5}
L(\sigma;k,X,n)=\sigma_1(k,X,n)
L_1(\sigma;k,X,n)+\sigma_2(k,X,n)L_2(\sigma;k,X,n)\bigr),
\end{equation}
where $L_1(\sigma;k,X,n)=L_2(\sigma;k,X,n)=0$ for $k+n=N$ and
\begin{gather}\label{b6}
\begin{aligned}
L_1(\sigma;k,X,n)&= \displaystyle{\int_0^C} M m \,\lambda(m| X,n)dm+L(\sigma;k+M,X,n),\\
L_2(\sigma;k,X,n)&= \displaystyle{\int_{-C}^0} \Big( |Mm|  \\
 &+\mE_Y L(\sigma;k,X+Y,n+M)\Big) \lambda(m |
X,n)dm
\end{aligned}
\end{gather}
for $k+n<N$. Here $L_\ell(\cdot)$ denotes cumulative expected
income on the residual control horizon of the length $(N-k-n)$ if
at first $M$ times the $\ell$-th action was chosen and then the
control was implemented according to strategy $\sigma$
($\ell=1,2$). Expected losses \eqref{a1} are equal to
$L(\sigma;0,0,0)$.

Let's transform equations~\eqref{b3}, \eqref{b4} and \eqref{b5},
\eqref{b6} to more convenient for calculations forms. Denote by
$f_D(X)= f_D(X|0)$ the density of Gaussian distribution with
variance $D$.

\begin{Theorem}\label{th21} Consider a dynamic programming
equation
\begin{equation}\label{b7}
\tilde{R}(k,X,n)=\min\bigl(\tilde{R}_1(k,X,n),\tilde{R}_2(k,X,n)\bigr),
\end{equation}
where $\tilde{R}_1(k,X,n)=\tilde{R}_2(k,X,n)=0$ for $k+n=N$ and
\begin{gather}\label{b8}
\begin{aligned}
\tilde{R}_1(k,X,n)&=M \tilde{g}_1(X,n)+\tilde{R}(k+M,X,n),\\
\tilde{R}_2(k,X,n)&=M \tilde{g}_2(X,n) \\
&+\int_{-\infty}^{+\infty} \tilde{R}(k,X+Y,n+M) h_{n,M}(MX-n Y)
\,dY
\end{aligned}
\end{gather}
for $0\le k+n<N$. Here
\begin{gather}\label{b9}
\begin{aligned}
\tilde{g}_1(X,n)&=\displaystyle{\int_0^C} m f_{n^*}(X|n m)\,\lambda(m)dm,\\
\tilde{g}_2(X,n)&=\displaystyle{\int_{-C}^0} |m| f_{n^*}(X|n
m)\,\lambda(m)dm.
\end{aligned}
\end{gather}
and
\begin{equation}\label{b10}
h_{0,M}(Y)=1,\qquad h_{n,M}(Y)= (n+M) f_{n M^* (n+M)} (Y),\quad
n>0.
\end{equation}

Then Bayesian risk~\eqref{a4} is equal to
\begin{equation}\label{b11}
R^B_N(\lambda)= \tilde{R}(0,0,0).
\end{equation}
\end{Theorem}

\begin{proof} Let's multiply \eqref{b4} by $P(X,n)$ defined in~\eqref{b2}.
Then the equation \eqref{b7}--\eqref{b8} holds for risks
$\tilde{R}(k,X,n)=R^B(k,X,n)P(X,n)$ with
\begin{gather*}
\begin{aligned}
h_{n,M}(MX-n Y)&= \frac{\int_{-C}^C f_{M^*}(Y|Mm) f_{n^*}(X|nm)
\,\lambda(m)dm}{\int_{-C}^C f_{n^*+M^*}(X+Y|(n+M) m)
\,\lambda(m)dm}\\ &= \displaystyle{\frac{f_{M^*}(Y|Mm) f_{n^*}(X|n
m)}{f_{n^*+M^*}(X+Y|(n+M) m)}}.
\end{aligned}
\end{gather*}
One can straightforwardly verify that this corresponds to
expressions in~\eqref{b10}. Clearly Bayesian risk~\eqref{a4} is
equal to~\eqref{b11}.
\end{proof}

Below we use the standard notation for convolution
$F(x)*G(x):=\int_{-\infty}^\infty F(x-y) G(y)\,dy$.

\begin{Theorem}\label{th22} Consider a dynamic programming
equation
\begin{equation}\label{b12}
R(k,X,n)=\min(R_1(k,X,n),R_2(k,X,n)),
\end{equation}
where $R_1(k,X,n)=R_2(k,X,n)=0$ for $k+n=N$ and
\begin{equation}\label{b13}
\begin{aligned}
R_1(k,X,n)&=M g_1(X,n)+R(k+M,X,n),\\ R_2(k,X,n)&=M
g_2(X,n)+R(k,X,n+M)* f_{M^*}\left(X\right)
\end{aligned}
\end{equation}
for $0\le k+n<N$. Here
\begin{equation}\label{b14}
\begin{aligned}
g_1(X,n) &=\int_0^C m \exp\left(D^{-1}(Xm-0.5 n m^2)\right)\,\lambda(m)dm,\\
g_2(X,n)&= \int_{-C}^0|m| \exp\left(D^{-1}(Xm-0.5 n
m^2)\right)\,\lambda(m)dm.
\end{aligned}
\end{equation}
Then Bayesian risk~\eqref{a4} is equal to
\begin{equation}\label{b15}
R^B_N(\lambda)= R(0,0,0).
\end{equation}
\end{Theorem}

\begin{proof} Let's put $\tilde{R}(k,X,n)=R(k,X,n)f_{n^*}(X)$. As
$\tilde{g}_\ell(X,n)\\ = g_\ell(X,n) f_{n^*}(X)$, $\ell=1,2$, the
first equation~\eqref{b13} follows from the first
equation~\eqref{b8}. The second equation~\eqref{b8} can be
transformed as follows
\begin{equation}\label{b16}
R_2(k,X,n)=M g_2(X,n)+\int_{-\infty}^\infty R(k,X+Y,n+M)
f_{M^*}(Y)\,dY.
\end{equation}
The validity of~\eqref{b16} for $n>0$ follows from the equality
$f_{n^*+M^*}(X+Y)h_{n,M}(MX-nY)=f_{n^*}(X)f_{M^*}(Y)$. For $n=0$
the equation~\eqref{b16} is verified by straightforward analysis.
The second equation~\eqref{b13} is equivalent to~\eqref{b16}
because of the evenness of $f_M(Y)$. Formula~\eqref{b15} follows
from~\eqref{b11}.
\end{proof}

\begin{Remark}\label{r3}
As $R(k,X,n)f_{n^*}(X)=R^B(k,X,n)P(X,n)$ then
\begin{gather*}
R(k,X,n)=R^B(k,X,n)\int_{-C}^C \exp\left(D^{-1}(Xm-0.5 n
m^2)\right)\,\lambda(m)dm.
\end{gather*}
\end{Remark}

Similar reasonings applied to equation \eqref{b5}--\eqref{b6}
result in the following theorem presented without proof.

\begin{Theorem}\label{th23} Given a strategy
$\sigma:\sigma_\ell(k,X,n)=\Pr(y_{k+n+1}=\ell|k,X,n)$,  let's
consider a dynamic programming equation
\begin{equation}\label{b17}
L(\sigma;k,X,n)=\sigma_1(k,X,n) L_1(\sigma;k,X,n)+\sigma_2(k,X,n)
L_2(\sigma;k,X,n),
\end{equation}
where $L_1(\sigma;k,X,n)=L_2(\sigma;k,X,n)=0$ for $k+n=N$ and
\begin{equation}\label{b18}
\begin{aligned}
L_1(\sigma;k,X,n)&=M g_1(X,n)+L(\sigma;k+M,X,n),\\
L_2(\sigma;k,X,n)&=M g_2(X,n)+L(\sigma;k,X,n+M)*
f_{M^*}\left(X\right)
\end{aligned}
\end{equation}
for $0\le k+n<N$.

Then expected losses~\eqref{a1} are equal to
\begin{equation}\label{b19}
L_N(\sigma;\lambda)= L(\sigma;0,0,0).
\end{equation}
\end{Theorem}

Considered problem in Bayesian setting for Bernoulli one-armed
bandit was previously investigated in~\cite{BJK,Cher}. In
particular, there were proved the thresholding property of
Bayesian strategy and its representation by two stages: at the
initial stage the second action only is applied and at the final
stage the first action only is applied. These two features of
Bayesian strategy are proved below in case of Gaussian
distributions.

\begin{Lemma}\label{lem1}
For all fixed $k,n$ the difference $\Delta
R(k,X,n)=R_1(k,X,n)-R_2(k,X,n)$ is monotonically increasing
function of $X$ and under condition~\eqref{b1} the following
equalities hold
\begin{equation}\label{b20}
\lim_{X\to -\infty}\Delta R(k,X,n)= -\infty,\qquad \lim_{X\to
+\infty}\Delta R(k,X,n)= +\infty.
\end{equation}
Therefore there exists $T(k,n)$ such that Bayesian strategy
chooses the first and the second actions according to conditions
$X<T(k,n)$ and $X>T(k,n)$ respectively (if $X=T(k,n)$ then actions
may be arbitrary chosen).
\end{Lemma}

\begin{Lemma}\label{lem2}
If $R_1(k,X,n)\le R_2(k,X,n)$ then $R_1(k,X,n)=(N-k-n)\times\,
g_1(X,n)$. Therefore, once being chosen, the first action will be
applied till the end of the control.
\end{Lemma}

Proofs of these lemmas are very close to those presented in
\cite{Koln15} and hence they are omitted here. From
Lemma~\ref{lem2} the following Theorem follows.

\begin{Theorem}\label{th24}
Let's denote statistics $(k,X,n)$, risks $R(k,X,n)$ and thresholds
$T(k,n)$ for $k=0$  by $(X,n)$, $R(X,n)$ and $T(n)$. It follows
from lemma~$\ref{lem2}$ that Bayesian strategy and Bayesian risk
can be searched for as a solution to the following recursive
equation
\begin{equation}\label{b21}
R(X,n)=\min(R_1(X,n),R_2(X,n)),
\end{equation}
where $R_1(X,n)=R_2(X,n)=0$ for $n=N$ and
\begin{equation}\label{b22}
\begin{aligned}
R_1(X,n)&=(N-n) g_1(X,n),\\ R_2(X,n)&=M g_2(X,n)+R(X,n+M)*
f_{M^*}(X),
\end{aligned}
\end{equation}
where $g_1(X,n),g_2(X,n)$ are defined in~\eqref{b14}. The first
action, once being chosen, will be applied till the end of the
control.
\end{Theorem}

\begin{Remark}
Batch strategy allows the following generalization. The data can
be partitioned into  $K$ batches of different sizes, so that the
data processing is implemented by batches of sizes
$M_1,\dots,M_{K}$, where $M_1+\dots+M_{K}=N$. Then equations
\eqref{b7}--\eqref{b8}, \eqref{b12}--\eqref{b13},
\eqref{b17}--\eqref{b18} and \eqref{b21}--\eqref{b22} hold true if
$M$ is replaced by $M_1$ at $n=0$ and by $M_i$ at
$n=M_1+\dots+M_{i-1}$ ($i
> 1$). It is expediently to choose batches of smaller sizes at the start of the control as the results
of Section~\ref{num} show. Some examples of usage of the batches
of different sizes are presented in~\cite{O17}.
\end{Remark}

The following theorem makes it possible to restrict consideration
with the variances $D = 1$.

\begin{Theorem}\label{th25}
Given some $k>0$, let the following transformations be made:
$\hD=k D$, $\hm=k^{1/2} m$, $\hX=k^{1/2} X$, $\hC=k^{1/2}C$,
$\hlambda(\hm)=k^{-1/2}\lambda(m)$, $\hM^*=kM^*$,
$\hsigma_\ell(k,\hX,n)=\sigma_\ell(k,X,n)$. Then corresponding
Bayesian risks and losses are related by equalities
\begin{gather}\label{b23}
 \hR^B_N(\hlambda)=k^{1/2} R^B_N(\lambda),
 \quad \hL_N(\hsigma,\hlambda)=k^{1/2} L_N(\sigma,\lambda).
\end{gather}
\end{Theorem}

\begin{proof}
Really, implementing the above transformations in
~\eqref{b12}--\eqref{b15} and in \eqref{b17}--\eqref{b19} one
obtains by induction that $\hR(k,\hX,n)=k^{1/2}R(k,X,n)$,
$\hL(k,\hX,n)=k^{1/2}L(k,X,n)$. Hence, the required equalities
follow from \eqref{b15}, \eqref{b19}.
\end{proof}

\begin{Corollary}\label{group}
Let's consider a treatment of $N=MK$ data by $K$ batches, each
containing $M$ data items, on the set of parameters
$\Theta_N=\{|m| \le C M^{-1/2}\}$ and one-by-one treatment of $K$
data items on the set of parameters $\Theta_K=\{|m| \le C \}$.
Then the following equality holds
\begin{equation}\label{b24}
N^{-1/2} R^M_{N(M)}(\Theta_N)= K^{-1/2 }R^M_{K(1)}(\Theta_K).
\end{equation}
Here in notations  $R^M_{N(M)}(\cdot)$ and $R^M_{K(1)}(\cdot)$ the
treatment by batches containing $M$ data items and one-by-one
treatment are explicitly indicated. Equality \eqref{b24} implies
that corresponding minimax risks depend only on the numbers of
processed batches.
\end{Corollary}
\begin{proof}
Let's put $m=\tm M^{-1/2}$. Then $\Theta_N=\{|\tm| \le C\}$, i.e.
this change of variables maps $\Theta_N$ on $\Theta_K$. Next, the
treatment of $K$ batches, each containing $M$ data items, is
equivalent to one-by-one treatment of $K$ data items with $\hD=M
D$, $\hm=M^{1/2} \tm$, the supports of the worst-case prior
distribution functions $\hlambda_0(\hm)$ and $\lambda_0(m)$ are
consistent with each other at that. By \eqref{b23} the validity of
equality $\hR^B_{N(M)}(\hlambda_0(\hm))=M^{1/2}
R^B_{K(1)}(\lambda_0(m))$ follows. This implies \eqref{b24}.
\end{proof}

According to \eqref{a3} the scaled Bayesian risk $N^{-1/2
}R^M_N(\Theta)$ is bounded from above. Hence, it follows from the
corollary~\ref{group} that batch data processing virtually does
not increase the minimax risk if the number of batches is large
enough. Let's suppose now that distributions of incomes are not
Gaussian. Nevertheless, according to the central limit theorem
distributions of cumulative incomes in large enough batches of
data are close to Gaussian. This implies that strategies of batch
data processing provide close values of the minimax risk for a
wide class of processes with equal mathematical expectations and
variances of one-step incomes, i.e. these strategies are
universal.

On the other hand, it follows from the corollary~\ref{group} that
batch data processing sets more restrictive requirements on the
set of parameters than one-by-one treatment. This is due to the
initial stage of control when possibly the worst action is applied
to the first batch of data and not to a single data item.

\section{Invariant Equations and Limiting Description}\label{lim}

Let's give an invariant notation of the equations describing batch
data processing with the control horizon equal to unit. Denote
$s=N^{-1} k$, $t=N^{-1} n$, $x= N^{-1/2} X$, $y= N^{-1/2} Y$,
$\eps=MN^{-1}$, $w=N^{1/2} m$, $\varrho(w)=N^{-1/2}\lambda(m)$,
$c=N^{1/2}C$, $r_{\ell,\eps}(s,x,t)= N^{-1/2} R_\ell(k,X,n)$,
$r_\eps(s,x,t)=N^{-1/2} R(k,X,n)$, $r_{\ell,\eps}(x,t)= N^{-1/2}
R_\ell(X,n)$, $r_\eps(x,t)=N^{-1/2} R(X,n)$. Consider the set of
parameters $\Theta_N=\{\theta:\:|w|\le c \}=\{\theta:\:|m|\le c
N^{-1/2}\}$ which describes the set of close distributions.

\begin{Theorem}\label{th31} Consider a dynamic programming
equation
\begin{equation}\label{c1}
r_\eps(s,x,t)=\min(r_{1,\eps}(s,x,t),r_{2,\eps}(s,x,t)),
\end{equation}
where $r_{1,\eps}(s,x,t)=r_{2,\eps}(s,x,t)=0$ for $s+t=1$ and
\begin{equation}\label{c2}
\begin{aligned}
r_{1,\eps}(s,x,t)&=\eps g_1(x,t)+r(s+\eps,x,t),\\
r_{2,\eps}(s,x,t)&=\eps g_2(x,t)+r_\eps(s,x,t+\eps)*
f_{\eps^*}\left(x\right)
\end{aligned}
\end{equation}
for $0\le k+n<N$. Here
\begin{equation}\label{c3}
\begin{aligned}
g_1(x,t) &=\int_0^c w \exp\left(D^{-1}(xw-0.5 t w^2)\right)\,\varrho(w)dw,\\
g_2(x,t)&= \int_{-c}^0|w| \exp\left(D^{-1}(xw-0.5 t
w^2)\right)\,\varrho(w)dw.
\end{aligned}
\end{equation}
Then Bayesian risk~\eqref{a4} is equal to
\begin{equation}\label{c4}
R^B_N(\lambda)= N^{1/2} r_\eps(0,0,0).
\end{equation}
\end{Theorem}

\begin{proof}
Theorem follows from theorem \ref{th22} by implementing the above
change of variables in \eqref{b12}--\eqref{b15}.
\end{proof}

\begin{Theorem}\label{th32} Given a strategy
$\sigma:\sigma_\ell(s,x,t)=\Pr(y_{s+t+\eps}=\ell|s,x,t)$,  let's
consider a dynamic programming equation
\begin{equation}\label{c5}
l_\eps(\sigma;s,x,t)=\sigma_1(s,x,t)
l_{1,\eps}(\sigma;s,x,t)+\sigma_2(s,x,t) l_{2,\eps}(\sigma;s,x,t),
\end{equation}
where $l_{1,\eps}(\sigma;s,x,t)=l_{2,\eps}(\sigma;s,x,t)=0$ for
$s+t=1$ and
\begin{equation}\label{c6}
\begin{aligned}
l_{1,\eps}(\sigma;s,x,t)&=\eps g_1(x,t)+l_\eps(\sigma;s+\eps,x,t),\\
l_{2,\eps}(\sigma;s,x,t)&=\eps g_2(x,t)+l_\eps(\sigma;s,x,t+\eps)*
f_{\eps^*}\left(x\right)
\end{aligned}
\end{equation}
for $0\le s+t<1$.

Then expected losses~\eqref{a1} are equal to
\begin{equation}\label{c7}
L_N(\sigma;\lambda)= N^{1/2}l_\eps(\sigma;0,0,0).
\end{equation}
\end{Theorem}

\begin{proof}
Theorem follows from theorem \ref{th23} by implementing the above
change of variables in \eqref{b17}--\eqref{b19}.
\end{proof}

\begin{Theorem}\label{th2}
Let's consider a dynamic programming equation
\begin{equation}\label{c8}
r_\eps(x,t)=\min(r_{1,\eps}(x,t),r_{2,\eps}(x,t)),
\end{equation}
where $r_{1,\eps}(x,1)=r_{2,\eps}(x,1)=0$ and
\begin{equation}\label{c9}
\begin{aligned}
r_{1,\eps}(x,t)&=(1-t)g_1(x,t),\\
r_{2,\eps}(x,t)&=\eps g_2(x,t)+r_\eps(x,t+\eps)*f_{\eps^*}(x)
\end{aligned}
\end{equation}
for $0\le t<1$, where $g_1(x,t)$, $g_2(x,t)$ are defined in
\eqref{c3}. Then Bayesian risk~\eqref{a4} is equal to
\begin{equation}\label{c10}
R^B_N(\lambda)= N^{1/2} r_\eps(0,0).
\end{equation}
\end{Theorem}

\begin{proof}
Theorem follows from theorem \ref{th24} by implementing the above
change of variables in \eqref{b21}--\eqref{b22}.
\end{proof}

In what follows we obtain some Lipschitz conditions for
$r_{\eps}(x,t)$. Let's make some preliminary remarks. For
$\theta=(0,wN^{-1/2})$, $n=tN^{1/2}$ denote by
\begin{gather}\label{c11}
l_\eps(\sigma,t,\rho)=N^{-1/2}\int_{-c}^c L_{N-n}(\sigma,\theta)
     \rho(w) dw, \quad r_\eps(t,\rho)=\inf_{\{\sigma\}} l_\eps(\sigma,t,\rho)
\end{gather}
scaled loss function and risk averaged with respect to the density
$\rho(w)$. We do not require that $\rho(w)$ is a probability
density, i.e. $\rho(w)\ge 0$ but it is possible that
\begin{gather*}
\int_{-c }^c \rho(w) dw \ne 1.
\end{gather*}

The following lemma holds.
\begin{Lemma}\label{le1} Given any densities
$\rho_1(w)$, $\rho_2(w)$, the inequalities hold
\begin{equation}
|r_\eps(t,\rho_1)-r_\eps(t,\rho_2)| \leq c (1-t) \int_{-c }^c
 |\rho_1(w)-\rho_2(w)|dw,\label{c12}
\end{equation}
\begin{equation}
\begin{aligned}
|r_\eps(t,&\rho_1)-r_\eps(t+\delta,\rho_2)|\\ & \leq  c
(1-t-\delta) \int_{-c }^c |\rho_1(w)-\rho_2(w)|dw + c\, \delta
\int_{-c }^c \rho_1(w)dw.\label{c13}
\end{aligned}
\end{equation}
\end{Lemma}

\begin{proof}
To check the validity of \eqref{c12} note that according to
\eqref{c11}
\begin{gather*}
\inf_{\{\sigma\}}l_\eps(\sigma,t,\rho_1)\le
\inf_{\{\sigma\}}l_\eps(\sigma,t,\rho_2)+ c(1-t)\int_{-c }^c
 |\rho_1(w)-\rho_2(w)|dw.
\end{gather*}
Since $\inf_{\{\sigma\}}l_\eps(\sigma,t,\rho)=r_\eps(t,\rho)$ and
inequality remains valid if one swaps $\rho_1$ and $\rho_2$ this
implies \eqref{c12}. Next, obviously
\begin{gather*}
r_\eps(t+\delta,\rho)\le r_\eps(t,\rho)\le
r_\eps(t+\delta,\rho)+c\delta \int_{-c }^c \rho(w)dw.
\end{gather*}
In view of \eqref{c12}, this implies \eqref{c13}.
\end{proof}

Let's denote $c'=D^{-1}c$.

\begin{Lemma}\label{lem3}
Functions $r_\eps(x,t),(r_\eps)'_x(x,t)$ are uniformly bounded.
For $t\ge \eps$ and arbitrary $x$ the following estimates hold
\begin{gather}
r_\eps(x,t)\le (1-t) \min\left(g_1(x,t),g_2(x,t)\right),\label{c14}\\
|(r_\eps)'_x (x,t)|\le c' r_\eps(x,t).\label{c15}
\end{gather}
For $t=0$, $x=0$ the following estimate holds
\begin{gather}
r_\eps(0,0)\le \eps g_2(0,0)+ (1-\eps)
f_{\eps^*}(x)*\min\left(g_1(x,\eps),g_2(x,\eps)\right),\label{c16}
\end{gather}
\end{Lemma}

\begin{proof}
Given statistics $(x,t)$, the one-step expected income for
choosing the $\ell$-th action is equal to $\eps g_\ell(x,t)$. The
estimate~\eqref{c14} is provided by the following strategy: at the
residual control horizon use only the action corresponding to the
larger value of $g_1(x,t)$, $g_2(x,t)$. And the estimate
~\eqref{c16} is ensured by the strategy: at the start of the
control use the second action and then at the residual control
horizon use only the action corresponding to the smaller value of
$g_1(x,\eps)$, $g_2(x,\eps)$.

Let's prove~\eqref{c15}. Note that
\begin{equation}\label{c17}
r_\eps(x,t)=\sigma_1(x,t)r_{1,\eps}(x,t)+\sigma_2(x,t)r_{2,\eps}(x,t),
\end{equation}
where $\sigma_\ell(x,t)$ is a probability to choose the $\ell$-th
action under given history $(x,t)$. In view of the thresholding
feature of the strategy $\sigma_\ell(x,t) \in\{0,1\}$ and its
switching at equality of $r_{1,\eps}(x,t)$, $r_{2,\eps}(x,t)$ we
obtain
\begin{equation}\label{c18}
|(r_\eps)'_x(x,t)|=\sigma_1(x,t)|(r_{1,\eps})'_x(x,t)|
+\sigma_2(x,t)|(r_{2,\eps})'_x(x,t)|.
\end{equation}
Since $|(g_\ell)'_x(x,t)|\le D^{-1}c g_\ell(x,t)=c' g_\ell(x,t)$,
$\ell=1,2$, it follows from~\eqref{c17}, \eqref{c18} that the
estimate~\eqref{c15} holds true for $t=1-\eps$.
Further~\eqref{c15} is proved by induction with the use of
obtained from~\eqref{c9} estimates
\begin{gather*}
|(r_{1,\eps})'_x(x,t)|\le (1-t)|(g_1)'_x(x,t)|\le c'
r_{1,\eps}(x,t),\\* |(r_{2,\eps})'_x(x,t)|\le\eps |(g_2)'_x (x,t)|
+|(r_\eps)'_x(x,t+\eps)|*f_{\eps^*}(x) \le c' r_{2,\eps}(x,t).
\end{gather*}
\end{proof}

\begin{Lemma}\label{lem4}
Assume that $\delta$ is a multiple of $\eps$.Then the following
estimate holds
\begin{equation}\label{c19}
|r_\eps(x,t)-r_\eps(x,t+\delta)|\le\delta  c \left(0.5 c c'
(1-t-\delta)+ 1 \right)\exp(c'|x|) .
\end{equation}
\end{Lemma}

\begin{proof}
Let's put $h(w,x,t)=\exp\left(D^{-1}(xw-0.5 t w^2)\right)$ and\\
$\rho_1(w)=h(w,x,t)\varrho(w)$,
$\rho_2(w)=h(w,x,t+\delta)\varrho(w)$. Note that $h(w,x,t)\le
\exp(c'|x|)$ and
\begin{equation}\label{c20}
|h'_t(w,x,t)|\le 0.5 D^{-1}c^2 \exp|D^{-1}cx|=0.5 c' c \exp(c'|x|)
\end{equation}
if $|w|\le c$. Using \eqref{c13}, \eqref{c20}  we obtain
\begin{gather*}
|r_{\eps}(x,t+\delta)- r_{\eps}(x,t)|\le \delta c (1-t-\delta)
\int_{-c}^c |h'_t(w,x,t)| d \varrho(w)\\ \qquad + \delta c
\int_{-c}^c h(w,x,t) d \varrho(w) \le \delta  c  \left(0.5 c c'
(1-t-\delta)+ 1 \right)\exp(c'|x|),
\end{gather*}
i.e. \eqref{c19} holds.
\end{proof}

Let's establish Lipschitz conditions with respect to $\varrho$,
for what we'll put $\varrho$ in notations $r_\eps(\varrho;x,t)$,
$g_\ell(\varrho;x,t)$. Define the distance between $\varrho_1$,
$\varrho_2$ by formula
\begin{gather*}
d(\varrho_1,\varrho_2)=\int_{-c}^c |\varrho_1(w)-\varrho_2(w)|dw.
\end{gather*}

\begin{Lemma}\label{lem5}
The following estimate holds
\begin{equation}\label{c21}
|r_\eps(\varrho_1;x,t)-r_\eps(\varrho_2;x,t)|\le (1-t)c_1(x)
\,d(\varrho_1,\varrho_2),
\end{equation}
with $c_1(x)=c \exp|c'x|$.
\end{Lemma}

\begin{proof}
Let's put $\rho_1(w)=h(w,x,t)\varrho_1(w)$,
$\rho_2(w)=h(w,x,t)\varrho_2(w)$. Using \eqref{c12} we obtain
\begin{gather*}
|r_\eps(\varrho_1;x,t)-r_\eps(\varrho_2;x,t)|\le c (1-t)
\int_{-c}^c h(w,x,t)| \varrho_1(w)- \varrho_2(w)|dw\\
\le  c (1-t) \exp(c'|x|) \int_{-c}^c | \varrho_1(w)-
\varrho_2(w)|dw,
\end{gather*}
i.e. \eqref{c21} holds.
\end{proof}

Let's consider now the approximation of the density $\varrho(w)$
by piece-wise constant density $\varrho^*(w)$. Denote by $\Delta
w=cK^{-1}$, $w_k=k\Delta w$, $k=-K,\ldots,K$. Let's put
\begin{gather*}
\varrho^*(w)=\varrho^*_k=\Delta w^{-1}\int_{w_k}^{w_{k+1}}
\varrho(w)dw \quad \mbox{for}\quad w\in (w_k,w_{k+1}).
\end{gather*}
Since $d(\varrho^*_1,\varrho^*_2)= \Delta w
\sum_{k=-K}^{K-1}|\varrho^*_{1,k}-\varrho^*_{2,k}|$, the densities
$\{\varrho^*\}$ are compact with respect to defined distance.

\begin{Lemma}\label{lem6}
The estimate holds
\begin{equation}\label{c22}
|r_\eps(\varrho^*;x,t)-r_\eps(\varrho;x,t)|\le c (1-t) c_2(x)
\Delta w,
\end{equation}
with $c_2(x)=\exp\left(c' |x|\right)D^{-1} (|x|+ c)$.
\end{Lemma}

\begin{proof}
Let's put $h(w,x,t)=\exp\left(D^{-1}(xw-0,5 t w^2)\right)$, then
\begin{gather*}
|h(w,x,t)'_w|\le \exp\left(D^{-1}c |x|\right)D^{-1}(|x|+ c)=
c_2(x).
\end{gather*}
In view of definition of $\varrho^*(w)$ and continuity of
$h(w,x,t)$, for some $w_k'\in (w_k,w_{k+1})$ the following
equality holds
$$
\int_{w_k}^{w_{k+1}} h(w,\cdot\,) \,\varrho(w)dw
=\int_{w_k}^{w_{k+1}} h(w_k',\cdot\,) \,\varrho^*(w)dw.
$$
Hence,
\begin{gather*}
\int_{w_k}^{w_{k+1}}|h(w,\cdot\,) \,\varrho(w)- h(w,\cdot\,)
\,\varrho^*(w)|dw\le c_2(x) \Delta w \int_{w_k}^{w_{k+1}}
\varrho^*(w)dw,
\end{gather*}
and consequently
\begin{gather*}
\int_{-c}^{c}|h(w,\cdot\,) \,\varrho(w)- h(w,\cdot\,)
\,\varrho^*(w)|dw\le c_2(x) \Delta w.
\end{gather*}
Then \eqref{c22} follows from \eqref{c12} with
$\rho_1(w)=h(w,x,t)\varrho(w)$,
\\ $\rho_1(w)=h(w,x,t)\varrho^*(w)$.

\end{proof}

Let's establish the existence of the limit of
$r_\eps(\varrho;x,t)$ as $\eps\to 0$.

\begin{Theorem}\label{th3}
For all $x,t$ for which the solution to equation~\eqref{c8},
\eqref{c9} is well defined there exist the limits
$r(\varrho;x,t)=\lim\limits_{\eps\to 0} r_\eps(\varrho;x,t)$ which
can be extended by continuity to all permissible $x,t$. These
limits are uniformly bounded, satisfy Lipschitz conditions with
respect to $x,t,\varrho$ and allow approximation by~$\varrho^*$
with constants presented in~\eqref{c14}, \eqref{c15}, \eqref{c19},
\eqref{c21} and \eqref{c22} respectively.

For the minimax risk on the set $\Theta_N=\{|m|\le cN^{-1/2}\}$
the estimate holds
\begin{equation}\label{c23}
\lim_{N\to\infty} N^{-1/2} R^M_N(\Theta_N)=\sup_\varrho
r(\varrho;0,0).
\end{equation}
\end{Theorem}

Proof of theorem is close to that presented in \cite{Koln15}.

A rigorous description of the limiting behavior of Bayesian risk
$r(x,t)$ turned out to be a difficult problem. So, let's present
nonrigourous reasonings and then supplement them
in~Section~\ref{num} by results of numerical experiments. Suppose
that $r_\eps(x,t)$ has partial derivatives of necessary orders and
show that the second equation~\eqref{c9} can be transformed to
\begin{equation}\label{c25}
r_{2,\eps}(x,t)=r_\eps(x,t+\eps)+ \eps \left( \frac{D}{2}\times
\frac{\partial^2 r_\eps(x,t+\eps)}{\partial
x^2}+g_2(x,t)\right)+o(\eps).
\end{equation}
For this purpose we represent $r_\eps(x-y,t+\eps)$ as Taylor
series:
\begin{equation}\label{c26}
r_\eps(x-y,\cdot\,)=r_\eps(x,\cdot\,)-y\times\frac{\partial
r_\eps(x,\cdot\,)}{\partial x}+\frac{y^2}{2}
\times\frac{\partial^2 r_\eps(x,\cdot\,)}{\partial x^2}+o(y^2).
\end{equation}
Let $A(\eps)=\{y:\:|y|\le\eps^{1/2-\alpha}\}$ with $0<\alpha<1/2$.
Taking into account that
$$
\begin{aligned}
\int_{\overline{A}(\eps)} f_\eps (y) \,dy&=o(\eps),\quad
&\int_{A(\eps)} f_\eps (y) \,dy&=1 +o(\eps),\\
\int_{A(\eps)} y f_\eps(y) \,dy&= o(\eps),\qquad &\int_{A(\eps)}
y^2 f_\eps(y) \,dy&= \eps+o(\eps),
\end{aligned}
$$
and substituting~\eqref{c26} in the second equation~\eqref{c9} we
obtain
$$
\begin{aligned}
r_{2,\eps}(x,t)&=\eps g_2(x,t)+ \int_{A(\eps)}
r_\eps(x-y,t+\eps) f_{\eps^*}(y) \,dy+ o(\eps)\\[-5pt] &=\eps
g_2(x,t)+ r_\eps(x,t+\eps)+\frac{D\eps}{2}\times \frac{\partial^2
r_\eps(x,t)}{\partial x^2}+o(\eps),
\end{aligned}
$$
i.e. \eqref{c25} holds true. The first equation~\eqref{c9} does
not vary at passage to the limit.

Continuing nonrigourous reasonings, note that if it is possible to
pass to the limit as $\eps \downarrow 0$ then two equations follow
from~\eqref{c9}, \eqref{c25} for $r=r(x,t)$:
$$
\begin{aligned}
(1-t) g_1(x,t)-r(x,t)&=0\quad \mbox{for}\quad (x,t)\in \Omega_1,\\
\frac{\partial r}{\partial t}+\frac{D}{2}\times \frac{\partial^2
r}{\partial x^2} +g_2(x,t)&=0\quad \mbox{for}\quad (x,t)\in
\Omega_2,
\end{aligned}
$$
where $\Omega_1,\Omega_2$ are domains in which the first and the
second actions are chosen respectively. Now let's recall that
these equations should be supplemented by equation~\eqref{c8}
which now can be written as
$$
\min_{\ell=1,2} (r_{\ell,\eps}(x,t)-r_\eps(x,t))=0,
$$
and then the limiting description of the function $r=r(x,t)$ takes
a form
\begin{equation}\label{c27}
\min \left((1-t) g_1(x,t)-r(x,t), \ \frac{\partial r}{\partial
t}+\frac{D}{2} \times \frac{\partial^2 r}{\partial x^2}
+g_2(x,t)\right)=0
\end{equation}
with initial and boundary conditions
\begin{equation}\label{c28}
r(x,1)=0,\qquad \lim_{x\to +\infty} r(x,t)= \lim_{x\to
-\infty}r(x,t)=0.
\end{equation}
Equation~\eqref{c27} describes function $r(x,t)$ and domains
$\Omega_1,\Omega_2$ together because the domain $\Omega_\ell$
corresponds to the minimum of the $\ell$-th entry in the left-hand
side of~\eqref{c27}. From~\eqref{c27} the difference equation
follows
\begin{equation}\label{c29}
r(x,t)=\min(r_1(x,t),r_2(x,t)),
\end{equation}
where
\begin{equation}\label{c30}
\begin{aligned}
r_1(x,t)&=(1-t)g_1(x,t),\\[-2pt] r_2(x,t)&=r(x,t+\Delta t)+ \Delta t \left(
\frac{D}{2}\times \Delta^2 r(x,t+\Delta t) +g_2(x,t)\right),
\end{aligned}
\end{equation}
initial and boundary conditions~\eqref{c28} are satisfied and
\begin{gather*}
\Delta^2 r(x,t) =\frac{r(x+\Delta x,t)-2r(x,t)+r(x-\Delta
x,t)}{\Delta x^2}.
\end{gather*}
The strategy providing a solution to
equation~\eqref{c29}--\eqref{c30} prescribes to choose the
$\ell$-th action if $r_\ell(x,t)$ has the minimum value.  The
first action, once being chosen, will be applied till the end of
the control.

The limiting expected losses $l(\sigma;0,0)$ are calculated as
follows
\begin{gather}\label{c31}
l(\sigma;x,t)=\sigma_1(x,t) l_1(\sigma;x,t) +\sigma_2(x,t)
l_2(\sigma;x,t),
\end{gather}
where
\begin{gather}\label{c32}
\begin{aligned}
l_1(\sigma;x,t)&=(1-t)g_1(x,t),\\[-2pt] l_2(\sigma;x,t)&=l(x,t+\Delta t)+ \Delta t
\left(\frac{D}{2}\times \Delta^2 l(x,t+\Delta t) +g_2(x,t)\right),
\end{aligned}
\end{gather}
with initial and boundary conditions
\begin{equation}\label{c33}
l(\sigma;x,1)=0,\qquad \lim_{x\to +\infty} l(\sigma;x,t)=
\lim_{x\to -\infty}l(\sigma;x,t)=0.
\end{equation}

\begin{Remark}
For the equation \eqref{c1}--\eqref{c2}, which is equivalent to
the equation \eqref{c8}--\eqref{c9}, corresponding partial
differential equation takes the form
\begin{equation}\label{c34}
\min \left(\frac{\partial r}{\partial s}+ g_1(x,t), \
\frac{\partial r}{\partial t}+\frac{D}{2}\times \frac{\partial^2
r}{\partial x^2} +g_2(x,t)\right)=0
\end{equation}
with initial and boundary conditions
\begin{equation}\label{c35}
r(s,x,t)|_{s+t=1}=0,\qquad \lim_{x\to +\infty} r(s,x,t)=
\lim_{x\to -\infty}r(s,x,t)=0.
\end{equation}
\end{Remark}

\section{Numerical Experiments}\label{num}

Minimax strategy and minimax risk were found by numerical methods
as Bayesian corresponding to the worst-case prior distribution
with the use of equation~\eqref{c29}, \eqref {c30} and
conditions~\eqref {c28} for $D=1$. It was assumed that the
worst-case prior distribution is concentrated at two points
$w=d_1$ and $w=-d_2$ with probabilities $\Pr(w=d_1)=\varrho$,
$\Pr(w=-d_2)=1-\varrho$. Clearly, $d_1,d_2,\varrho$ should
correspond to the maximum of $r(\varrho;0,0)$. So, they were
determined as $d_1\approx 1.65$, $d_2\approx 2.52$,
$\varrho\approx 0.38$ with $\max r(\varrho;0,0)\approx 0.37$. We
applied $\Delta t=5000^{-1}$, $\Delta x=0.0143$ for calculations
($\Delta t$, $\Delta x$ must satisfy the inequality $\Delta t/
\Delta x^2<1$).

\begin{figure}
\includegraphics[height=8cm]{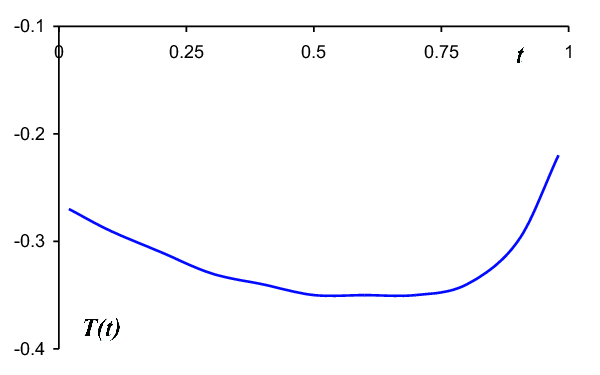}
\caption[]{Thresholds of the minimax strategy.} \label{fig1}
\end{figure}

Determined strategy has a thresholding property. It prescribes to
apply the second action if $x>T(t)$ and to switch to the first
action till the end of the control if  $x<T(t)$, where function
$T(t)$ describes threshold values presented on Figure~\ref{fig1}.

\begin{figure}
\includegraphics[height=8cm]{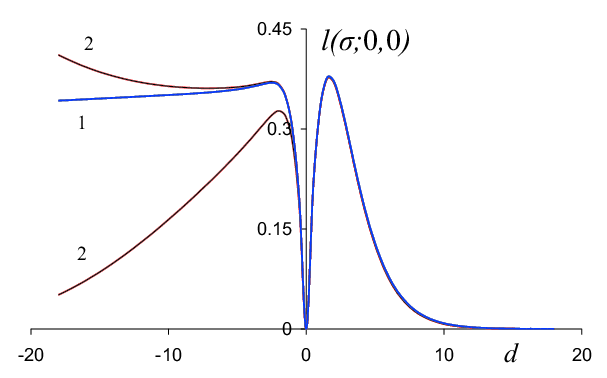}
\caption[]{The impact of the initial stage on expected losses.}
\label{fig2}
\end{figure}

\begin{figure}
\includegraphics[height=8cm]{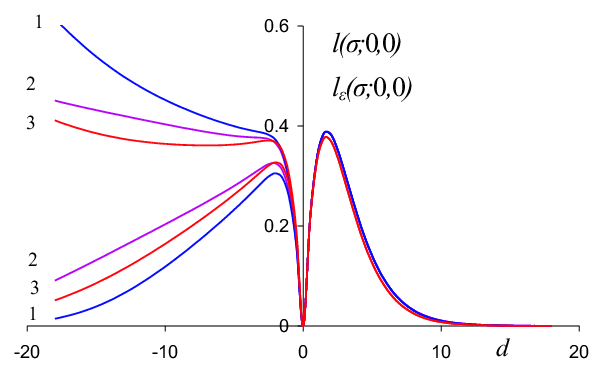}
\caption[]{The impact of the size of the batch on expected
losses.} \label{fig3}
\end{figure}

Then for determined strategy expected losses were calculated. The
case $d>0$ corresponds to $d_1=d$, $\varrho=1$ and the case $d<0$
corresponds to $d_2=-d$, $\varrho=0$. Results are presented on
Figure~\ref{fig2}. Curve~1 describes expected losses
$l(\sigma;0,0)$ determined according to \eqref{c31}--\eqref{c33}.
Expected losses described by the curve~1 have two maxima at
$d\approx 1.65$ and at $d\approx-2.52$ which are approximately
equal to $0.37$ and this confirms the made assumption of the
worst-case prior distribution and that determined Bayesian
strategy and risk are minimax ones. Note that according to the
implemented substitution of variables these maxima correspond to
the values of mathematical expectations $m=d_1 N^{-1/2}$ and
$m=-d_2 N^{-1/2}$. This means that maximal losses are attained for
close values of parameters if $N$ is large enough.

\begin{figure}
\includegraphics[height=8cm]{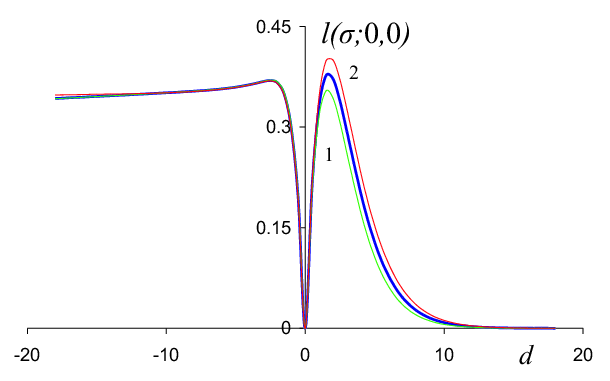}
\caption[]{To the robustness of the strategy–1.} \label{fig4}
\end{figure}

Curves~2 on Figure~\ref{fig2} are obtained for the strategy that
at the start of the control uses the first action for $\eps=0.02$
fraction of the control horizon. If $d>0$ then these losses are
almost the same as $l(\sigma;0,0)$ presented by curve~1. If $d<0$
then upper curve 2 describes expected losses
$|d|\eps+l(\sigma;x,\eps)*f_{\eps^*}(x)$ including those at the
initial stage. Lower curve~2 describes expected losses
$l(\sigma;x,\eps)*f_{\eps^*}(x)$ without those at the initial
stage. Since $l(\sigma;x,\eps)*f_{\eps^*}(x)\to 0$  with growing
$|d|$ then expected losses
$|d|\eps+l(\sigma;x,\eps)*f_{\eps^*}(x)$ exceed the Bayesian risk
$r(\varrho;0,0)$ at approximately $|d|\eps>r(\varrho;0,0)$. In the
domain $|d|\le \eps^{-1}r(\varrho;0,0)\approx 0.37 \times 50=18.5$
the strategy is the minimax one if $d<0$.

On Figure~\ref{fig3} curves~3 are the same as curves~2 on
Figure~\ref{fig2}. Curves~1 and 2 present batch data processing
and are calculated by equations~\eqref{c8}--\eqref{c10}. Upper
curves~1 and 2 describe expected losses
$|d|\eps+l_\eps(\sigma;x,\eps)*f_{\eps^*}(x)$ including those at
the initial stage. Lower curves~1 and 2 describe expected losses
$l(\sigma;x,\eps)*f_{\eps^*}(x)$ without those at the initial
stage. Curves~1 and 2 are calculated for $\eps=30^{-1}$ and
$\eps=50^{-1}$ respectively. If $d>0$ then all the curves are
close to each other. Note that $l(\sigma;0,0)$ and
$l_\eps(\sigma;0,0)$ describe scaled losses of batch processing
implemented in the infinitely many
 and in $\eps^{-1}$ stages respectively.  For calculations we
chosen step of integration $0.002$, integral limits from $-1.6$ to
$+1.6$.

\begin{figure}
\includegraphics[height=8cm]{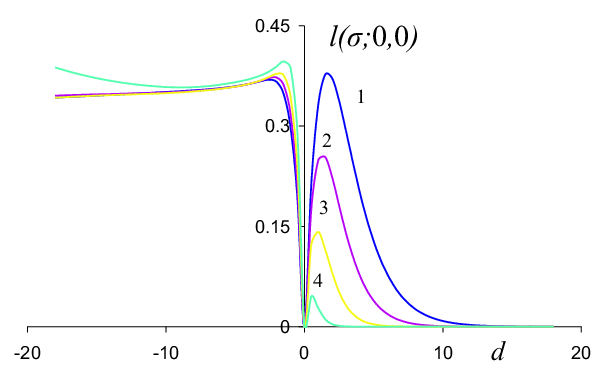}
\caption[]{To the robustness of the strategy-2.} \label{fig5}
\end{figure}

Figures~\ref{fig4} and~\ref{fig5} demonstrate the impact of the
variance $D$ on the losses. On Figure~\ref{fig4} the thick curve
presents $l(\sigma;0,0)$ corresponding to $D=1$. Thin curves~1 and
2 present expected losses corresponding to determined minimax
strategy for $D=1$ if actually $D=0.95$ and $D=1.05$ respectively.
One can see that all the curves are close to each other. Given a
large number of data, this means that one can estimate the
variance at the initial stage of control when the second action
only is applied and then use the estimate for calculating the
minimax strategy.

On Figure~\ref{fig5} the curve~1 presents $l(\sigma;0,0)$
calculated for $D=1$ if actually $D=1$. And curves~2, 3 and 4
present $l(\sigma;0,0)$ corresponding to determined minimax
strategy if actually $D=0.75,\ 0.5,\ 0.25$ respectively. This
means that the minimax strategy calculated for some $D_0$ remains
to be close to  minimax if $0.25 D_0\le D \le D_0$.

\begin{figure}
\includegraphics[height=8cm]{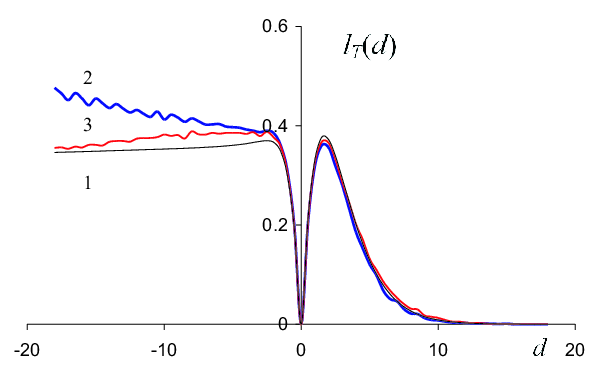}
\caption[]{Monte-Carlo simulation results.} \label{fig6}
\end{figure}

Finally, on Figure~\ref{fig6} Monte-Carlo simulation results are
presented for batch processing of $T=5000$ Bernoulli incomes
implemented in $N=50$ stages by batches of $M=100$ items with the
use of determined minimax strategy. Probabilities of successful
and unsuccessful processing by the first and the second methods
were equal to $p_1=p=0.5$ and $p_2=p+d(D/T)^{1/2}$ respectively,
where $D=p(1-p)=0.25$. The control of partitioned data was
implemented by the method described in~Section~1 but $\{\xi_n\}$
are scaled by the factor $D^{1/2}$ so that the variances of
$\{D^{-1/2}\xi_n\}$ are approximately equal to unit.

On Figure~\ref{fig6} the curve~1 presents $l(\sigma;0,0)$. The
curve~2 describes Monte-Carlo simulation results for the scaled
expected losses
\begin{gather*}
l_T(d)=(DT)^{-1/2}\left(T (p \vee p_2)-\mE_{\sigma,\theta}
\left(\sum_{t=1}^T \zeta_t\right)\right).
\end{gather*}
One can see that the curve~2 grows approximately linearly if $d$
is negative and its absolute value grows. This is due to the
losses for the first batch of 100 data processing. To avoid this,
the data on the first two stages were processed by eight batches
each containing 25 items and then the size of the batch was again
equal to 100. This processing corresponds to the curve~3 and there
is no the growth of losses there.

\section{Comparison with One-by-One Processing. Bernoulli Case}\label{asymp}

Results of the previous section imply that batch data processing
almost does not increase the minimax risk of the Gaussian
one-armed bandit if the number of batches is large enough.
Moreover, according to the central limit theorem these results
imply that batch data processing provides the close value of the
minimax risk for a wide class of one-armed bandits with the same
mathematical expectations and variances of one-step incomes.

However, there is a question of control performance in case of
one-by-one data processing. Is it possible to diminish the value
of the minimax risk in this case? In this section, we consider the
asymptotic (as $N \to \infty$) description of Bayesian risk of
Bernoulli one-armed bandit and show that it obeys the same partial
differential equation as the Bayesian risk of Gaussian one-armed
bandit. The reasonings below are not quite rigorous because we do
not have estimates on the smoothness of Bayesian risk.

Let's consider a Bernoulli one-armed bandit which is described by
a controlled random process $\xi_n$, $n=1,2,\dots,N$. Incomes
$\{\xi_n\}$ depend only on the currently chosen action $y_n$ as
follows
\begin{gather*}
    \Pr(\xi_n=1|y_n=\ell)=p_\ell, \quad \Pr(\xi_n=0|y_n=\ell)=q_\ell,
\end{gather*}
$\ell=1,2$. A Bernoulli one-armed bandit is described by a
parameter $\theta=(p,p_2)$ where $0<p<1$ is assumed to be known.
The set of parameters is as follows $\Theta=\{\theta=(p,p_2):
|p_2-p|\le C)\}$.

A control strategy $\sigma$ at the point of time $k+n+1$ assigns,
in general, a random selection of actions depending on the current
history which is described by triplet $(k,X,n)$, where $k,n$ are
current cumulative numbers of the first and the second actions
applications, $X$ is current cumulative income for the use of the
second action. In what follows we assume that at the start of the
control the strategy uses the second action $n_0$ times. This
almost does not influence the Bayesian risk if $n_0 \ll N$ and can
only a little increase it. The loss function is defined as
\begin{gather*}
L_N(\sigma,\theta)=N (p \vee
p_2)-\mE_{\sigma,\theta}\Biggl(\,\sum_{n=1}^N \xi_n\Biggr)
\end{gather*}
Given a prior probability distribution density $\lambda(p_2)$,
Bayesian risk is defined as follows
\begin{equation}\label{d1}
R^B_N(\lambda)= \min_{\{\sigma\}} \int_\Theta L_N(\sigma,\theta)
\,\lambda(p_2) dp_2.
\end{equation}

Let's obtain a standard dynamic programming equation for
calculating the Bayesian risk \eqref{d1}. The posterior
distribution density corresponding to the history of the process
$(k,X,n)$ is the following
\begin{gather}\label{d2}
\begin{aligned}
\lambda(p_2|X,n)&=\displaystyle{\frac{B(X,n|p_2)\ \lambda(p_2)} {P(X,n)}},\\
\mbox{where}\quad &P(X,n)= \displaystyle{\int_{p-C}^{p+C}}
B(X,n|p_2)\,\lambda(p_2)dp_2.
\end{aligned}
\end{gather}
Here
\begin{gather*}
  B(X,n|p)=\binom{n}X p^X q^{n-X}, \quad \mbox{where }\ q=1-p.
\end{gather*}
If additionally it is assumed that $B(X,n|p_2)=1$ as $n=0$ then
$\lambda(p_2|0,0)=\lambda(p_2)$.

Denote by $R^B(k,X,n)$ Bayesian risk on the last $(N-k-n)$ steps
calculated with respect to the posterior distribution
$\lambda(p_2| X,n)$. To find the Bayesian risk~\eqref{d1} one has
to solve the following standard recursive equation
\begin{equation}\label{d3}
R^B(k,X,n)=\min\bigl(R^B_1(k,X,n),R^B_2(k,X,n)\bigr),
\end{equation}
where $R^B_1(k,X,n)=R^B_2(k,X,n)=0$ for $k+n=N$ and
\begin{gather}\label{d4}
\begin{aligned}
R^B_1(k,X,n)&= \displaystyle{\int_p^{p+C}}  (p_2-p) \,\lambda(p_2| X,n)dp_2+R^B(k+1,X,n),\\
R^B_2(k,X,n)&= \displaystyle{\int_{p-C}^p} \Big( |p_2-p|\\
 &+\mE_Y R^B(k,X+Y,n+1)\Big) \lambda(p_2|X,n)dp_2
\end{aligned}
\end{gather}
for $k+n<N$. Here $\mE_Y$ denotes the mathematical expectation:
$$
\mE_Y R(Y)= \sum_{Y\in \{0,+1\}} R(Y) p_2^Y q_2^{1-Y}=q_2 R(0)+p_2
R(1).
$$
In equations above, $R^B_\ell(\cdot)$ means cumulative expected
income on the residual control horizon of the length $(N-k-n)$ if
at first the $\ell$-th action was chosen and then the control was
optimally implemented ($\ell=1,2$). Bayesian strategy prescribes
always to choose the action corresponding to the smaller value of
$R^B_1(\cdot),R^B_2(\cdot)$; the choice may be arbitrary in case
of their equality. In view of the above assumption on the
strategy, the Bayesian risk \eqref{d1} is equal to
\begin{gather}\label{d5}
    R^B(\lambda)=n_0 \int_{p-C}^p |p_2-p| \lambda(p_2)dp_2
+\int_{-\infty}^\infty R^B(0,X,n_0)P(X,n_0) dX .
\end{gather}

Denote $\tR_1(k,X,n)=R^B_1(k,X,n)P(X,n)$,
$\tR_2(k,X,n)\\=R^B_2(k,X,n)P(X,n)$, $\tR(k,X,n)=R^B(k,X,n)P(X,n)$
with $P(X,n)$ defined in~\eqref{d2}. Using \eqref{d2}--\eqref{d4}
one obtains
\begin{equation}\label{d6}
\tilde{R}(k,X,n)=\min\bigl(\tilde{R}_1(k,X,n),\tilde{R}_2(k,X,n)\bigr),
\end{equation}
where $\tilde{R}_1(k,X,n)=\tilde{R}_2(k,X,n)=0$ for $k+n=N$ and
\begin{equation}\label{d7}
\begin{aligned}
\tilde{R}_1(k,X,n)&=\tilde{g}_1(X,n)+\tilde{R}(k+1,X,n),\\
\tilde{R}_2(k,X,n)&=\tilde{g}_2(X,n) +\tilde{R}(k,X,n+1)\times
h(X,n,0)\\ &+ \tilde{R}(k,X+1,n+1)\times h(X,n,1)
\end{aligned}
\end{equation}
for $0\le k+n<N$. Here
\begin{equation}
\begin{aligned}\label{d8}
&\tg_1(X,n)=\int_p^{p+C} (p_2-p) B(X,n| p_2) \lambda(p_2)dp_2,\\
&\tg_2(X,n)=\int_{p-C}^p |p_2-p| B(X,n|p_2)\,\lambda(p_2)dp_2,
\end{aligned}
\end{equation}
and
\begin{equation}\label{d9}
h(X,n,0)=\frac{n+1-X}{n+1},\qquad h(X,n,1)=\frac{X+1}{n+1},\quad
\mbox{if } n>0.
\end{equation}
Let's check the validity of~\eqref{d9}. Analysis of the second
equations~~\eqref{d4}, \eqref{d7} gives that
$$
h(X,n,Y)= \frac{\int_{-C}^C  B(X,n|p_2)p_2^Y q_2^{1-Y}
\,\lambda(p_2)dp_2}{\int_{-C}^C B(X+Y,n+1|p_2)
\,\lambda(p_2)dp_2}=
\frac{\displaystyle{\binom{n}{X}}}{\displaystyle{\binom{n+1}{X+Y}}}.
$$
One can directly check that this corresponds to expressions
in~\eqref{d9}. Clearly, Bayesian risk~\eqref{d5} is equal to
\begin{equation}\label{d10}
R^B_N(\lambda)=n_0 \int_{p-C}^p |p_2-p| \lambda(p_2)dp_2
+\int_{-\infty}^\infty \tilde{R}(0,X,n_0) dX.
\end{equation}

Now let's assume that $N$, $n_0$ are large enough and $n_0 \ll N$.
Consider the following change of variables: $p_1=p$,
$p_2=p+wN^{-1/2}$. $C=c N^{-1/2}$,
$\lambda(p_2)=N^{1/2}\varrho(w)$. $X=np+x N^{1/2}$. $t=nN^{-1}$.
$s=kN^{-1}$, $\eps=N^{-1}$, $\delta=N^{-1/2}$. If $n>n_0$ and
$n_0$ is large enough then according to the central limit theorem
 $B(X,n|p_2)=f_{n^*}(X|np_2)(1+o(1))$. Here
\begin{gather*}
f_{n^*}(X|np_2)=(2\pi
nD)^{-1/2}\exp\left(-\frac{(X-np_2)^2}{2nD}\right)\\=N^{-1/2}(2\pi
tD)^{-1/2}\exp\left(-\frac{(x-wt)^2}{2Dt}\right)=N^{-1/2}f_{t^*}(x|wt).
\end{gather*}
This and \eqref{d8} imply that
\begin{equation}\label{d11}
\tilde{g}_1(X,n)=\eps \tilde{g}_1(x,t)+o(\eps),\quad
\tilde{g}_2(X,n)=\eps \tilde{g}_2(x,t)+o(\eps).
\end{equation}
with
\begin{equation}\label{d12}
\begin{aligned}
\tilde{g}_1(x,t)&= \int_0^c w f_{t^*}(x|wt) \varrho(w)dw,\\
\tilde{g}_2(x,t)&= \int_{-c}^0 |w| f_{t^*}(x|wt) \varrho(w)dw.
\end{aligned}
\end{equation}

Let's put  $\tilde{R}(k,X,n)= \tilde{r}(s,x,t)$. In view of
\eqref{d11} the first equation~\eqref{d7} takes the form
\begin{equation}\label{d13}
\tilde{r}_1(s,x,t)=\eps \tilde{g}_1(x,t)+
\tilde{r}(s+\eps,x,t)+o(\eps).
\end{equation}

To write in new variables the second equation~\eqref{d7} note that
$x'$ and $x$ corresponding to points of time $n+1$ and $n$ are
related as
\begin{equation*}
x'-x=\frac{X-(n+1)p}{N^{1/2}}-\frac{X-np}{N^{1/2}}=-p\delta
\end{equation*}
and hence $\tilde{R}(k,X,n+1)= \tilde{r}(s,x-p\delta,t+\eps)$.
Similarly,
\begin{equation*}
x'-x=\frac{X+1-(n+1)p}{N^{1/2}}-\frac{X-np}{N^{1/2}}=q\delta
\end{equation*}
and hence $\tilde{R}(k,X+1,n+1)=\tilde{r}(s,x+q\delta,t+\eps)$.
Additionally,
\begin{gather*}
\frac{n+1-X}{n+1}=q+\frac{p}{n+1}-\frac{xN^{1/2}}{n+1}=q+p
t^{-1}\eps -x t^{-1} \delta +o(\eps),\\[5pt]
\frac{X+1}{n+1}=p+\frac{q}{n+1}+\frac{xN^{1/2}}{n+1}=p+q
t^{-1}\eps +x t^{-1} \delta +o(\eps).
\end{gather*}
So, the second equation~\eqref{d7} takes the form
\begin{gather}\label{d14}
\begin{aligned}
\tilde{r}_2(s,x,t)&=\eps \tilde{g}_2(x,t) +
\tilde{r}(s,x-p\delta,t+\eps)\times (q+p t^{-1}\eps-x
t^{-1}\delta)\\ &+  \tilde{r}(s,x+q\delta,t+\eps)\times (p+q
t^{-1}\eps+x t^{-1}\delta)+o(\eps).
\end{aligned}
\end{gather}

The Bayesian risk \eqref{d3} is equal to
\begin{gather}\label{d15}
    R^B(\lambda)=N^{1/2}\left(\eps_0 \int_{-c}^0 |w|
    \varrho(w)dw
+\int_{-\infty}^\infty \tr(0,x,\eps_0) dx +o(\eps)\right).
\end{gather}

If $\tr(s,x,t)$ is smooth enough and has partial derivatives of
the necessary orders then equations~\eqref{d13} and \eqref{d14}
take the form
\begin{gather}\label{d16}
\begin{aligned}
\tr_1(s,x,t)&=\eps \tg_1(x,t)+ \tr(s+\eps,x,t)+o(\eps),\\
\qquad\tr_2(s,x,t)&= \eps \tg_2(x,t)+
\eps\tr(s,x,t+\eps)\times(1+t^{-1})\\ &+\tr'_x(s,x,t+\eps) x
t^{-1} +0.5\eps \tr''_{xx}(s,x,t+\eps)D +o(\eps).
\end{aligned}
\end{gather}
with $D=pq$. Note that \eqref{d16} is not rigorously derived
because we do not have estimates on the smoothness of
$\tr(s,x,t)$. Equations \eqref{d16} must be complemented by
equation \eqref{d6} written in new variables $s,x,t$. Now we
present~\eqref{d6} as
\begin{equation}\label{d17}
\min_{\ell=1,2}(\tr_\ell(s,x,t)-\tr(s,x,t)=0
\end{equation}

From~\eqref{d16}, \eqref{d17} one obtains in the limiting case (as
$\eps \to 0$) the second order partial differential equation
\begin{gather}\label{d18}
\begin{aligned}
\min \Big(\tg_1(x,t)&+ \tr'_s(s,x,t), \quad
\tg_2(x,t)+ \tr'_t(s,x,t)\\
 &+t^{-1}\tr(s,x,t)+ x t^{-1} \tr'_x(s,x,t)  +0.5 D
\tr''_{xx}(s,x,t)\Big)=0.
\end{aligned}
\end{gather}

The following theorem holds.

\begin{Theorem}\label{th41}
Let $\tr(s,x,t)$ satisfy the second order partial differential
equation~\eqref{d18}. Then it can be expressed as
\begin{gather}\label{d19}
\tilde{r}(s,x,t)=r(s,x,t)f_{t^*}(x),
\end{gather}
where $r(s,x,t)$ satisfies the second order partial differential
equation
\begin{gather}\label{d20}
\min \left(\frac{\partial r}{\partial s}+ g_1(x,t), \frac{\partial
r}{\partial t}+ \frac{D}{2} \times\frac{\partial^2 r}{\partial
x^2}+ g_2(x,t) \right)=0,
\end{gather}
and $g_1(x,t)$, $g_2(x,t)$ are defined in~\eqref{c3}. Bayesian
risk is expressed as
\begin{equation}\label{d21}
R^B_N(\lambda)=N^{1/2} \left(\eps_0 \int_{-c}^0 |w| \varrho(w)dw
+\int_{-\infty}^\infty r(0,x,\eps_0) f_{\eps_0^*}(x) dx\right).
\end{equation}

\end{Theorem}

\begin{proof}
Let's omit the dependence of $f$, $r$, $g_1$, $g_2$ on $s,x,t$, so
that $f'_t$ denotes now partial derivative by $t$. Then
\begin{gather}\label{d22}
\tg_1+ \tr'_s= (g_1+r'_s)f.
\end{gather}
Next,
\begin{gather}\label{d23}
\begin{aligned}
\tg_2+ \tr'_t
 &+t^{-1}\tr+ x t^{-1} \tr'_x  +0.5 D\tr''_{xx}(s,x,t)=g_2  f + (r'_t f +r f'_t)\\ &+ t^{-1} r f+ x t^{-1} (r'_x f +r f'_x)
+0.5 D (r''_{xx}f +2 r'_x f'_x +r f''_{xx}).
\end{aligned}
\end{gather}
Using the equalities
\begin{gather*}
f'_t=\left(-\frac{1}{2t}+\frac{x^2}{2Dt^2}\right)f, \quad
f'_x=-\frac{x}{Dt}f, \quad f''_{xx}=2 D^{-1} f'_t,
\end{gather*}
one derives that \eqref{d23} is equal to
\begin{gather}\label{d24}
\begin{aligned}
(g_2+r'_t&+t^{-1} r + x t^{-1} r'_x + 0.5 D r''_{xx})f+2 r f'_t +
x t^{-1} r f'_x+D r'_x f'_x\\  = \Big( g_2&+r'_t+t^{-1} r + x
t^{-1} r'_x + 0.5 D r''_{xx}+2 r
\Big(-\displaystyle{\frac{1}{2t}}+\displaystyle{\frac{x^2}{2Dt^2}}\Big)
\\ & - r \displaystyle{\frac{x^2}{Dt^2}}
-\displaystyle{\frac{x}{t}} r'_x\Big) f= \Big( g_2+r'_t+ 0.5 D
r''_{xx} \Big) f.
\end{aligned}
\end{gather}
From \eqref{d22}, \eqref{d24} one derives \eqref{d20}. Formula
\eqref{d21} follows from \eqref{d15} and \eqref{d19}.
\end{proof}

One can see that equations \eqref{c34} and \eqref{d20} which
describe Bayesian risks for batch data processing and for
one-by-one processing are the same. Corresponding Bayesian risks
given by formulas~\eqref{c4} and \eqref{d21} are close to each
other if $\eps_0$ is small enough. This means that one-by-one data
processing cannot diminish Bayesian risk ensured by the batch data
processing if the number of batches is large enough. According to
the main theorem of the theory of games it cannot diminish the
minimax risk, too.

\section{Discussion}\label{disc}

Gaussian one-armed bandit naturally arises when optimization of
the batch data processing is considered if there are two
alternative processing methods available with a priori unknown
efficiency of the second method. In this case, the same methods
are applied to all the data in the same batches and then
cumulative incomes are used for the control. It turned out that
batch data processing almost does not increase the minimax risk if
the number of batches is large enough. For example, the scaled
minimax risk is only about 3\% higher its limiting value in case
of processing the data by 50 batches.

Since distributions of cumulative incomes in sufficiently large
batches of data are close to Gaussian for a wide class of random
processes with equal mathematical expectations and variances of
one-step incomes, the proposed strategies are universal. Moreover,
it seems highly likely that even one-by-one optimal data
processing does not allow to diminish the limiting scaled value of
the minimax risk. In the article, this assumption is confirmed for
Bernoulli two-armed bandit.

Proposed strategies demonstrate fine robustness properties. It
turned out that expected losses deviate just a little if the
variance is assigned with a significant error up to 5\%. This
makes it possible to estimate the variance at initial stage of the
control and then use the obtained estimate for determining the
strategy. In addition, there are no high requirements to closeness
of distributions of cumulative incomes in batches to Gaussian. For
example, Monte-Carlo simulations in Section~\ref{num} were
implemented for batches of 100 or 25 items of data. These
properties also imply universality of proposed strategies.
%%%%%%%%%%%%%%%%%%%%%%%%%%%%%%%%%%%%%%%%%%


\vspace{6pt}

\begin{thebibliography}{9}

\bibitem{BF}
\textsc{Berry, D. A.} and \textsc{Fristedt, B.} (1985).
\textit{Bandit Problems: Sequential Allocation of Experiments},
Chapman~\&\ Hall, London.

\bibitem{PS}
\textsc{Presman, E. L.} and \textsc{Sonin, I.M.} (1990).
\textit{Sequential Control with Incomplete Information: Bayesian
Approach}, Academic Press, New York.

\bibitem{Tsetlin}
\textsc{Tsetlin, M. L.} (1973). \textit{Automaton Theory and
Modeling of Biological Systems}, Academic Press, New York.

\bibitem{Sragovich}
\textsc{Sragovich, V. G.} (2006). \textit{Mathematical Theory of
Adaptive Control}, World Sci., Singapore.


\bibitem{Robbins}
\textsc{Robbins, H.} (1952). Some Aspects of the Sequential Design
of Experiments. \textit{Bull.\ Amer.\ Math.\ Soc.} \textbf{58}
527--535.

\bibitem {FZ}
\textsc{Fabius J.} and  \textsc{van~Zwet W. R.} (1970). Some
Remarks on the Two-Armed Bandit. \textit{Ann.\ Math.\ Statist}.
\textbf{41} 1906--1916.

\bibitem{Vogel}
\textsc{Vogel, W.} (1960). An Asymptotic Minimax Theorem for the
Two Armed Bandit Problem. \textit{Ann.\ Math.\ Stat.}  \textbf{31}
444--451.

\bibitem{Nazin}
\textsc{Juditsky, A.}, \textsc{Nazin, A. V.}, \textsc{Tsybakov, A.
B.} and \textsc{Vayatis, N.} (2008) Gap-Free Bounds for Stochastic
Multi-Armed Bandit. In \textit{Proc.\ 17th IFAC World Congr.,
Seoul, Korea, July~6--11, 2008}. 11560--11563. Available at
http://www.ifac-papersonline.net/Detailed/37644.html.

\bibitem{LLRS}
\textsc{Lai, T.L.}, \textsc{Levin, B.}, \textsc{Robbins, H.} and
\textsc{Siegmund, D.} (1980). Sequential Medical Trials.
\textit{Proc.\ Natl.\ Acad.\ Sci.\ USA.} \textbf{77} 3135--3138.

\bibitem{Koln11}
\textsc{Kolnogorov, A. V.} (2011) Finding Minimax Strategy and
Minimax Risk in a Random Environment (the Two-Armed Bandit
Problem). \textit{Automation and Remote Control}. \textbf{72}
1017--1027.

\bibitem{BJK}
\textsc{Bradt, R. N.}, \textsc{Johnson, S. M.} and \textsc{Karlin,
S.} (1956) On Sequential Designs for Maximizing the Sum of\/ $n$
Observations. \textit{Ann.\ Math.\ Statist.} \textbf{27}
1060--1074.

\bibitem{Cher}
\textsc{Chernoff, H.} and \textsc{Ray, S. N.} (1965) A Bayes
Sequential Sampling Inspection Plan. \textit{Ann.\ Math.\
Statist.} \textbf{36} 1387--1407.

\bibitem{Koln15}
\textsc{Kolnogorov, A. V.} (2015) One-Armed Bandit Problem for
Parallel Data Processing Systems. \textit{Problems of Information
Transmission}  \textbf{51} 177--191.

\bibitem{Koln18}
\textsc{Kolnogorov, A. V.} (2018) Gaussian Two-Armed Bandit and
Optimization of Batch Data Processing. \textit{Problems of
Information Transmission} \textbf{54} 84--100.

\bibitem {O17}
\textsc{Oleynikov, A. O.} (2013) Numerical Optimization of
Parallel Processing in a Stationary Environment. \textit{Trans.\
Karelian Res.\ Centre Russ.\ Acad.\ Sci.} \textbf{1} 73--78 (in
russian).

\end{thebibliography}
\end{document}